\documentclass[11pt]{article}

\usepackage[utf8]{inputenc}
\usepackage[a4paper,margin=2.6cm]{geometry}
\usepackage{amsmath,amssymb,amsfonts,amsthm}
\usepackage{mathrsfs}
\usepackage{mathtools}
\usepackage{booktabs}
\usepackage{array}
\usepackage[hidelinks]{hyperref}

\theoremstyle{plain}
\newtheorem{theorem}{Theorem}[section]
\newtheorem{proposition}[theorem]{Proposition}
\newtheorem{lemma}[theorem]{Lemma}
\newtheorem{corollary}[theorem]{Corollary}
\theoremstyle{remark}
\newtheorem{remark}[theorem]{Remark}

\numberwithin{equation}{section}
\DeclareTextSymbolDefault{\textsection}{OMS}

\title{Nine-Colored Partition Components Modulo 19}
\author{G\"oksal Bilgici\\[4pt]
\small Department of Elementary Mathematics Education, Kastamonu University,\\
\small Kastamonu, T\"urkiye\\
\small \texttt{gbilgici@kastamonu.edu.tr}\quad
ORCID: \href{https://orcid.org/0000-0001-9964-5578}{0000-0001-9964-5578}}
\date{}

\begin{document}

\maketitle

\begin{abstract}
Let $p_{m}(n)$ be the number of partitions of $n$ into parts of $m$ colors, and let $y=q^{19}$. Starting from Winquist's identity, as Hirschhorn did for the modulus $11$, we express each of the nineteen components of $(q;q)_{\infty}^{10}$ modulo $19$ as a sum of two quartic products of the theta functions $(y^{a};y^{19})_{\infty}(y^{19-a};y^{19})_{\infty}$, and deduce all nineteen generating functions $\sum_{n} p_{9}(19n+r)y^{n}$ modulo $19$ in closed form. The vanishing class $p_{9}(19n+17)$ is the one known previously. For seventeen colors, the residue class singled out by Kiming and Olsson carries no congruence, and we show that its generating function is $y\,f(y^{19})^{2}/f(y)$ modulo $19$. Along the way, we recall, with proofs, such exact relationships known for all these congruences. O'Brien's cubic relation is a Weierstrass equation for $X_{0}(19)$ in Kolberg's coordinates; we prove the existence of a doubling automorphism of the theta functions, fixing the cubic, which acts on both dissections by explicit scalars. All proofs are by hand.
\end{abstract}



\section{Introduction}\label{sec:intro}

For an integer $m$, let $p_{m}(n)$ be defined by
\begin{equation}
\sum_{n\ge0} p_{m}(n)q^{n} = (q;q)_{\infty}^{-m},
\end{equation}
so that $p_{1} = p$ is the partition function and, for $m\ge1$, $p_{m}(n)$ counts partitions of $n$ into parts of $m$ colors. Following Kiming and Olsson \cite{KimingOlsson1992}, one asks for which primes $\ell$, numbers of colors $m$, and residues $a$ the congruence $p_{m}(\ell n+a) \equiv 0 \pmod{\ell}$ holds for all $n$. Since $(q;q)_{\infty}^{\ell} \equiv (q^{\ell};q^{\ell})_{\infty} \pmod{\ell}$, the numbers $p_{\ell-k}(n)$ are governed modulo $\ell$ by the $\ell$-dissection of $(q;q)_{\infty}^{k}$, and the question is about which components of that dissection vanish. 

Kiming and Olsson showed that, apart from the cases $k=1$ and $k=3$ coming from Euler's and Jacobi's identities, a congruence for $p_{\ell-k}$ can occur in only one single class $24a \equiv \ell-k$, and they listed all of them for $\ell\le19$ \cite[Theorems 1 and 4]{KimingOlsson1992}. Boylan \cite{Boylan2004} explained most of these congruences through the lacunary powers of Dedekind's $\eta$-function classified by Serre \cite{Serre1985}, and Locus and Wagner \cite{LocusWagner2017} gave a related treatment. At $\ell=19$ the list is
\begin{equation}\label{eq:KO19}
p_{7}(19n+9) \equiv 0, \qquad p_{9}(19n+17) \equiv 0, \qquad p_{13}(19n+14) \equiv 0 \pmod{19}.
\end{equation}

These describe single residue class statements only. In this paper, we find complete components using elementary methods from classical $q$-series identities. Put $y=q^{19}$ and $P_{a}=(y^{a};y^{19})_{\infty}(y^{19-a};y^{19})_{\infty}$.

\medskip\noindent\textbf{Nine colors.} 
Winquist's identity gives a double series representation \cite{Winquist1969} for $(q;q)_{\infty}^{10}$. We slice it modulo $19$ by pairing its indices, and verify that the component in the class $r$ is $P_{0}^{2}B_{r}$, where $B_{r}$ is a sum of two quartic products $P_{A} P_{B} P_{A+B} P_{B-A}$ with explicit coefficients, the pair $\{A, B\}$ falling into the class $r \equiv 11(A^{2}+B^{2}) - 2$. A count of the solutions to $A^{2}+B^{2} \equiv c$ shows that every class receives exactly two pairs, except for class $17$, which gets none (Theorem~\ref{thm:tenth} and Corollary~\ref{cor:two-terms}). This yields all nineteen of the components of the nine-colored partition function:
\[
\sum_{n\ge0}p_{9}(19n+r)\,y^{n}\equiv
\frac{P_{0}\,B_{r}}{P_{1}P_{2}\cdots P_{9}}\pmod{19}\qquad
(0\le r\le18),
\]
which is expressed as the sum of two quotients of products (Theorem~\ref{thm:colored} and Table~\ref{tab:Psi}). The second congruence in \eqref{eq:KO19} is the case $r=17$. Hirschhorn's method derives from Winquist's identity the $11$-dissection of $(q;q)_{\infty}^{10}$ and with it all eleven components of $p(n)$ modulo $11$, recovering the Atkin--Swinnerton-Dyer congruences for that modulus \cite{Hirschhorn2000}, \cite[\S13.4]{Hirschhorn2017}.

The counterpart at $19$ for nine colors, given in Theorem~\ref{thm:colored}, requires two products in place of one. What changes is the shape of the answer: at $11$ every nonzero component of $(q;q)_{\infty}^{10}$ is a single product, while at $19$ none is, and the count above shows why (Remark~\ref{rmk:contrast}).

\medskip\noindent\textbf{Seventeen colors.} 
Here the relevant power is $(q;q)_{\infty}^{2}$, and the class singled out by Kiming and Olsson is that of $11$. The case $k=2$ does not appear among the congruences described in \cite[Theorem~2.3]{Boylan2004}, and indeed $19$ is not exceptional for seventeen colors \cite[Theorem~1.3]{Boylan2004}. Instead, we show what that class is (Corollary~\ref{cor:seventeen}):
\[
\sum_{n\ge0}p_{17}(19n+11)\,y^{n}\equiv\frac{y\,f(y^{19})^{2}}{f(y)}
\pmod{19},
\]
where $f(y)=(y;y)_{\infty}$; specifically $p_{17}(19n+11)\equiv p(n-1)$ for $1\le n\le19$. At each prime $p \equiv 3 \pmod 4$ an analogous formula is valid, and there are two congruences of the same kind governed by $\sigma_{3}$ and $\sigma_{7}$. Closed forms of this kind for whole components go back to Newman \cite{Newman1957b}, in other residue classes.

\medskip\noindent\textbf{What is known.} 
The second and third congruences of \eqref{eq:KO19}, for nine and thirteen colors, the family $p_{p-10}(pn+5(p^{2}-1)/12)\equiv0\pmod p$ for primes $p\equiv3\pmod4$ that contains the second, and the exact relations $c(p^{2}s+\delta)=p^{h}c(s)$ behind it and behind the seventeen-color formula, are all known: they are Hecke relations at inert primes for forms with complex multiplication by $\mathbb{Q}(i)$ \cite{Serre1985,Boylan2004}. For $(q;q)_{\infty}^{10}$ and $(q;q)_{\infty}^{2}$, the exact relations go back to Newman \cite{Newman1953,Newman1955,Newman1957a} and were proved in general, by the argument used here, by Cooper, Hirschhorn, and Lewis \cite[Theorem~1]{CHL2000}.

We include the proofs (Sections~\ref{sec:family} and~\ref{sec:seventeen}), in our notation, to keep the paper self-contained, and state at each point whose they are. The congruence $p_{7}(19n+9)\equiv0$ is the one that Boylan calls superexceptional; we do not prove it, but we reformulate it in Section~\ref{sec:questions}.

\medskip\noindent\textbf{The exceptional component of $p(n)$.} 
The component $\Phi_{4}=\sum p(19n+4)y^{n}$ was studied by Kolberg \cite{Kolberg1959} and determined by O'Brien \cite{OBrien1965}. O'Brien's work contains the cubic $S=P_{4}P_{6}P_{9}+yP_{1}P_{7}P_{8}-y^{2}P_{2}P_{3}P_{5}$, the doubling symmetry that produces it, the congruence $\Phi_{4}\equiv5S/(P_{0}\,\mathcal{D}^{3})$ with $\mathcal{D}=P_{1}\cdots P_{9}$, and an identity over $\mathbb{Z}$ which reads $S^{3}=8y\mathcal{D}^{2}S+\mathcal{D}^{5}+19y^{3}\mathcal{D}$. We observe that in terms of Kolberg's basis functions this identity is a Weierstrass equation for the modular curve $X_{0}(19)$ (Remark~\ref{rmk:weierstrass}). Modulo $19$ the curve becomes nodal, and the cubic expresses $\mathcal{D}/P_{0}^{4}$ as $S$ times the explicit series $(S^{2}-8y\mathcal{D}^{2})/f(y)^{4}$; since $S$ is a unit in $\mathbb{F}_{19}[[y]]$, the content of this lies in the cofactor.

Our contribution here is to show that this structure is compatible with the two dissections. In the polynomial ring $\mathcal{R}=\mathbb{F}_{19}[Z_{1},\dots,Z_{9}]$, where $Z_{a}$ stands for $q^{t_{a}}P_{a}$, the two dissections become polynomials $B$ and $C$, and the weight-$11$ part of $BC^{3}$ has value $q^{30}\mathcal{D}/P_{0}^{4}$ (Proposition~\ref{thm:scalar}). O'Brien's symmetry becomes an automorphism $\mathcal{T}$ of $\mathcal{R}$ that preserves Weierstrass's three-term theta relations, multiplies $B$, $C$, and that weight-$11$ part by explicit scalars, and fixes the cubic (Lemmas~\ref{lem:T-ideal} and~\ref{lem:T-eigen}, Proposition~\ref{prop:S}). These facts are consistent with a factorization of that weight-$11$ part through the cubic in the polynomial ring, but do not prove one; that is the first of the questions in Section~\ref{sec:questions}.

\medskip\noindent\textbf{Method.} 
Every proof in the paper is by hand, from the classical identities of Euler, Jacobi, Winquist, Weierstrass, and Ramanujan, and from the results of Kolberg and O'Brien quoted with their equation numbers. No computer is required: the two tables are evaluations of direct formulas, which can be done manually. This dissection of Winquist's series follows Hirschhorn's treatment of the modulus $11$ \cite{Hirschhorn2000}, and was used for the modulus $17$ in \cite{Bilgici17}.

\medskip\noindent\textbf{Outline.} 
Section~\ref{sec:prelim} fixes notation and dissects $(q;q)_{\infty}^{3}$ and $(q;q)_{\infty}^{10}$ modulo $19$. Section~\ref{sec:nine} gives the nine-colored components, and Section~\ref{sec:family} the family containing $p_{9}(19n+17)$. Section~\ref{sec:seventeen} treats seventeen, sixteen, and thirteen colors. Section~\ref{sec:bridge} combines the two dissections, and Section~\ref{sec:exceptional} treats the exceptional component. Section~\ref{sec:questions} lists open questions.

\section{Notation and the two dissections}\label{sec:prelim}

Throughout, $f(q)=\prod_{n\ge1}(1-q^{n})$ and $y=q^{19}$. Unless otherwise noted, a congruence between two power series means congruence of their coefficients modulo $19$. The numbers $p_{m}(n)$ are defined by $\sum_{n \ge 0} p_{m}(n) q^{n}=f(q)^{-m}$, so that $p_{1} = p$ is the partition function and, for $m\ge1$, $p_{m}(n)$ counts partitions of $n$ into parts of $m$ colors.

For every $F\in\mathbb{Z}[[q]]$, there exists a unique expansion of the form
\begin{equation}\label{eq:dissection}
F(q)=\sum_{r=0}^{18}q^{r}F_{r}(y),\qquad F_{r}\in\mathbb{Z}[[y]],
\end{equation}
and we refer to $F_{r}$ as the \emph{component} of $F$ at level $r$. Given $G\in\mathbb{Z}[[y]]$, we have $(GF)_{r}=G\,F_{r}$; this fact is used repeatedly without comment.

\subsection{Products and theta series}

Put
\[
P_{0}=(y^{19};y^{19})_{\infty},\qquad
P_{a}=(y^{a};y^{19})_{\infty}(y^{19-a};y^{19})_{\infty}\quad(1\le a\le18),
\]
so that $P_{19-a}=P_{a}$, and let $\mathcal{D}=P_{1}P_{2}\cdots P_{9}$.
Factoring $f(y) = (y;y)_{\infty}$ according to the residue of the exponents modulo $19$ gives
\begin{equation}\label{eq:f-of-y}
f(y)=P_{0}\,\mathcal{D}.
\end{equation}

For an integer $x$, let $\rho(x)\in\{0,1,\dots,9\}$ be the representative of $\pm x$ modulo $19$.

The triple product identity of Jacobi, namely
\begin{equation}\label{eq:JTP}
\sum_{m\in\mathbb{Z}}(-1)^{m}x^{m(m-1)/2}w^{m}
=(w;x)_{\infty}(x/w;x)_{\infty}(x;x)_{\infty},
\end{equation}
gives, with $x=y^{19}$ and $w=y^{a}$,
\begin{equation}\label{eq:P-as-theta}
\sum_{m\in\mathbb{Z}}(-1)^{m}y^{19m(m-1)/2+am}=P_{0}P_{a}
\qquad(1\le a\le18).
\end{equation}

We also need theta series not of one conductor, but of three times the conductor. For every integer $h$, put
\begin{equation}\label{eq:J-def}
J_{h}=J_{h}(y)=\sum_{m\in\mathbb{Z}}(-1)^{m}y^{57m(m-1)/2+hm},
\end{equation}
which is a Laurent series in $y$. Replacing $m$ by $-m$ in the series for $J_{57-h}$ and $1-m$ in $J_{-h}$ gives
\begin{equation}\label{eq:J-sym}
J_{57-h}=J_{h},\qquad J_{-h}=-y^{-h}J_{h},\qquad J_{h+57}=-y^{-h}J_{h},
\end{equation}
the third relation following as a consequence of the first two. Finally, for integers $u$ and $v$, put
\begin{equation}\label{eq:t-and-d}
t_{u}=\frac{(9-u)(10-u)}{2},\qquad d(v)=\frac{(3-v)(10-3v)}{2},
\end{equation}
both of which are integers, and observe that
\begin{equation}\label{eq:t-d-shift}
t_{-u}-t_{u}=19u,\qquad d(-v)-d(v)=19v.
\end{equation}

\subsection{The cube}

\begin{proposition}[Jacobi's identity, dissected]\label{prop:jacobi}
\[
f(q)^{3}\equiv P_{0}\sum_{a=1}^{9}(-1)^{a}\,2a\,q^{t_{a}}P_{a}\pmod{19}.
\]
Here
\[
(t_{1},\dots,t_{9})=(36,28,21,15,10,6,3,1,0)\equiv(17,9,2,15,10,6,3,1,0)
\pmod{19},
\]
and these residues are distinct. Hence nine components of $f(q)^{3}$ are
single products, and the other ten vanish modulo $19$.
\end{proposition}

\begin{proof}
Jacobi's identity is
$f(q)^{3}=\sum_{n\ge0}(-1)^{n}(2n+1)q^{n(n+1)/2}$. The summand is unchanged
by $n\mapsto-1-n$, so $2f(q)^{3}$ is the same sum over all $n\in\mathbb{Z}$.
Write $n=19m+9-u$ with $m\in\mathbb{Z}$ and $-9\le u\le9$. Then
$2n+1\equiv-2u$, $(-1)^{n}=(-1)^{m+u+1}$, and
\[
\frac{n(n+1)}{2}=t_{u}+19\Bigl(\frac{19m(m-1)}{2}+(19-u)m\Bigr).
\]
Modulo $19$, the terms with a given $u$ therefore sum to
\[
(-1)^{u}\,2u\,q^{t_{u}}\sum_{m\in\mathbb{Z}}(-1)^{m}y^{19m(m-1)/2+(19-u)m}.
\]
For $u=0$ this is $0$. For $1\le u\le9$ the inner sum is $P_{0}P_{19-u}=
P_{0}P_{u}$ by \eqref{eq:P-as-theta}. The involution $n\mapsto-1-n$
exchanges the classes $u$ and $-u$, so they contribute equally, and the
factor $2$ cancels.
\end{proof}

\subsection{The tenth power}

Two identities of Winquist are used. The first is his product identity
\cite{Winquist1969}; see also \cite{Hirschhorn2000,Hirschhorn2017}. For
nonzero $a$ and $b$,
\begin{multline}\label{eq:winquist-2var}
(a,q/a,b,q/b,ab,q/(ab),a/b,bq/a;q)_{\infty}\,(q;q)_{\infty}^{2}\\
=\sum_{m,n\in\mathbb{Z}}(-1)^{m+n}q^{\frac{3m(m+1)}{2}+\frac{n(3n+1)}{2}}
\Bigl(a^{-3m}b^{-3n}-a^{-3m}b^{3n+1}\\
-a^{-3n+1}b^{-3m-1}+a^{3n+2}b^{-3m-1}\Bigr),
\end{multline}
where $(a_{1},\dots,a_{k};q)_{\infty}=\prod_{i}(a_{i};q)_{\infty}$. The
second is the series for the tenth power that he deduced from it:
\begin{equation}\label{eq:winquist-f10}
48\,f(q)^{10}=\sum_{a,b\in\mathbb{Z}}(-1)^{a+b}\,XY\,(X^{2}-Y^{2})\,
q^{\frac{3a(a+1)}{2}+\frac{b(3b+1)}{2}},
\qquad X=6a+3,\ \ Y=6b+1 .
\end{equation}
The exponent $N$ of $q$ in \eqref{eq:winquist-f10} satisfies
\begin{equation}\label{eq:norm-10}
24N+10=X^{2}+Y^{2}.
\end{equation}

For $1\le A<B\le9$ define the quartic product
\[
G_{A,B}=P_{A}\,P_{B}\,P_{\rho(A+B)}\,P_{B-A},
\]
and for every integer $v$ put $K_{v}=J_{19+3v}-y^{v}J_{19-3v}$.

\begin{lemma}\label{lem:winquist-special}
For $1\le A<B\le9$,
\[
P_{0}^{2}\,G_{A,B}=J_{3B}K_{A}-y^{B-A}J_{3A}K_{B}.
\]
\end{lemma}

\begin{proof}
In \eqref{eq:winquist-2var} replace $q$ by $y^{19}$ and put $a=y^{A}$,
$b=y^{B}$. On the right each of the four double sums is a product of two
single sums, and by \eqref{eq:J-def} and \eqref{eq:J-sym}
\begin{gather*}
\sum_{m}(-1)^{m}y^{\frac{57m(m+1)}{2}-3Am}=J_{57-3A}=J_{3A},\\
\sum_{n}(-1)^{n}y^{\frac{19n(3n+1)}{2}\mp3Bn}=J_{38\mp3B}=J_{19\pm3B}.
\end{gather*}
The right side becomes
\begin{multline*}
J_{3A}J_{19+3B}-y^{B}J_{3A}J_{19-3B}-y^{A-B}J_{3B}J_{19+3A}
+y^{2A-B}J_{3B}J_{19-3A}\\
=J_{3A}K_{B}-y^{A-B}J_{3B}K_{A}.
\end{multline*}
On the left, the first six factors give $P_{A}$, $P_{B}$ and
$P_{\rho(A+B)}$, because $(y^{c},y^{19-c};y^{19})_{\infty}=P_{c}=
P_{\rho(c)}$ for $1\le c\le18$. For the last two, with $d=B-A>0$,
\[
(y^{-d};y^{19})_{\infty}(y^{19+d};y^{19})_{\infty}
=\frac{1-y^{-d}}{1-y^{d}}\,(y^{19-d};y^{19})_{\infty}(y^{d};y^{19})_{\infty}
=-y^{-d}P_{d}.
\]
So the left side is $-y^{A-B}P_{0}^{2}G_{A,B}$. Multiplying both sides by
$-y^{B-A}$ gives the lemma.
\end{proof}

The next theorem is proved following Hirschhorn's method for dissecting $f(q)^{10}$ into $11$ pieces \cite{Hirschhorn2000}, \cite[\S13.4]{Hirschhorn2017}. Each of the single summations in \eqref{eq:winquist-f10} is dissected individually, then the resulting terms are paired up using Lemma~\ref{lem:winquist-special}. When dissecting with respect to $11$, all the pieces of a component are received by one pair. When dissecting with respect to $19$, each non-zero component receives two of these pairs (Corollary~\ref{cor:two-terms}).

\begin{theorem}[The tenth power modulo $19$]\label{thm:tenth}
For $1\le A<B\le9$ put
\[
\varepsilon(A,B)=3t_{B}+d(A)=\frac{3A^{2}+3B^{2}-19A-57B+300}{2},
\]
and let $r(A,B)\in\{0,\dots,18\}$ and $\kappa(A,B)\ge0$ be defined by
$\varepsilon(A,B)=19\,\kappa(A,B)+r(A,B)$. Then
\begin{equation}\label{eq:r-of-pair}
r(A,B)\equiv11\,(A^{2}+B^{2})-2\pmod{19},
\end{equation}
and
\[
f(q)^{10}\equiv P_{0}^{2}\sum_{r=0}^{18}q^{r}B_{r}(y)\pmod{19},
\]
where
\[
B_{r}=-3\!\!\sum_{\substack{1\le A<B\le9\\ r(A,B)=r}}\!\!
(-1)^{A+B}AB\,(B^{2}-A^{2})\;y^{\kappa(A,B)}\,G_{A,B}.
\]
\end{theorem}

\begin{proof}
For $k\in\{1,3\}$ let
\[
\mathcal{A}_{k}=\sum_{a\in\mathbb{Z}}(-1)^{a}(6a+3)^{k}q^{3a(a+1)/2},
\qquad
\mathcal{B}_{k}=\sum_{b\in\mathbb{Z}}(-1)^{b}(6b+1)^{k}q^{b(3b+1)/2},
\]
so that \eqref{eq:winquist-f10} reads
$48f(q)^{10}=\mathcal{A}_{3}\mathcal{B}_{1}-\mathcal{A}_{1}\mathcal{B}_{3}$.

\emph{Step 1: the series $\mathcal{A}_{k}$.} Write $a=19m+9-u$ with
$m\in\mathbb{Z}$ and $-9\le u\le9$. Then $6a+3\equiv-6u$,
$(-1)^{a}=(-1)^{m+u+1}$, and
\[
\frac{3a(a+1)}{2}=3t_{u}+19\Bigl(\frac{57m(m-1)}{2}+(57-3u)m\Bigr).
\]
Since $k$ is odd, the terms with a given $u$ sum, modulo $19$, to
\[
(-1)^{u}6^{k}u^{k}\,q^{3t_{u}}J_{57-3u}=(-1)^{u}6^{k}u^{k}\,q^{3t_{u}}J_{3u}.
\]
By \eqref{eq:J-sym} and \eqref{eq:t-d-shift},
$q^{3t_{-u}}J_{-3u}=q^{3t_{u}}y^{3u}\cdot(-y^{-3u})J_{3u}
=-q^{3t_{u}}J_{3u}$, while $(-1)^{-u}(-u)^{k}=-(-1)^{u}u^{k}$; so the classes
$u$ and $-u$ contribute equally, and
\begin{equation}\label{eq:A-k}
\mathcal{A}_{k}\equiv2\cdot6^{k}\sum_{u=1}^{9}(-1)^{u}u^{k}\,q^{3t_{u}}J_{3u}.
\end{equation}

\emph{Step 2: the series $\mathcal{B}_{k}$.} Write $b=19m+3-v$ with
$m\in\mathbb{Z}$ and $-9\le v\le9$. Then $6b+1\equiv-6v$,
$(-1)^{b}=(-1)^{m+v+1}$, and
\[
\frac{b(3b+1)}{2}=d(v)+19\Bigl(\frac{57m(m-1)}{2}+(38-3v)m\Bigr).
\]
The terms with a given $v$ sum to
\[
(-1)^{v}6^{k}v^{k}\,q^{d(v)}J_{38-3v}=(-1)^{v}6^{k}v^{k}\,q^{d(v)}J_{19+3v}.
\]
By \eqref{eq:t-d-shift} the class $-v$ contributes
$-(-1)^{v}6^{k}v^{k}\,q^{d(v)}y^{v}J_{19-3v}$, so
\begin{equation}\label{eq:B-k}
\mathcal{B}_{k}\equiv6^{k}\sum_{v=1}^{9}(-1)^{v}v^{k}\,q^{d(v)}K_{v}.
\end{equation}

\emph{Step 3: pairing.} From \eqref{eq:A-k} and \eqref{eq:B-k},
\[
48f(q)^{10}\equiv2\cdot6^{4}\sum_{u,v=1}^{9}(-1)^{u+v}uv\,(u^{2}-v^{2})\,
q^{3t_{u}+d(v)}J_{3u}K_{v}.
\]
The terms with $u=v$ vanish. For $A<B$ the two terms $(u,v)=(B,A)$ and
$(u,v)=(A,B)$ have opposite coefficients, and
$\bigl(3t_{A}+d(B)\bigr)-\bigl(3t_{B}+d(A)\bigr)=19(B-A)$. Together they
give
\[
(-1)^{A+B}AB\,(B^{2}-A^{2})\,q^{\varepsilon(A,B)}
\bigl(J_{3B}K_{A}-y^{B-A}J_{3A}K_{B}\bigr),
\]
and the bracket is $P_{0}^{2}G_{A,B}$ by
Lemma~\ref{lem:winquist-special}. Finally $2\cdot6^{4}/48=54\equiv-3$,
and $q^{\varepsilon}=q^{r}y^{\kappa}$. The congruence \eqref{eq:r-of-pair}
follows from $2\varepsilon(A,B)\equiv3(A^{2}+B^{2})+15$ and
$2^{-1}\equiv10$.
\end{proof}

\begin{corollary}\label{cor:two-terms}
$B_{17}=0$, and for every $r\ne17$ the sum defining $B_{r}$ has exactly two
terms. The eighteen nonzero $B_{r}$ are listed in Table~\ref{tab:B}.
\end{corollary}

\begin{proof}
Using \eqref{eq:r-of-pair} and $11 \cdot 7 \equiv 1$, we can determine the pairs contributing to $B_{r}$ as those with $A^{2}+B^{2} \equiv c$, where $c = 7(r+2)$. Since $19 \equiv 3 \pmod 4$, a sum of two nonzero squares is nonzero modulo $19$, thus no pair has $c = 0$, that is, $r = 17$. Fix $c \neq 0$, and let $\chi$ be the Legendre symbol modulo $19$. The equation $x^{2}+y^{2}=c$ has $19 - \chi(-1) = 20$ solutions in $\mathbb{F}_{19}^{2}$ \cite[Ch.~8]{IrelandRosen}. 

Here, $2(1 + \chi(c))$ solutions have $xy = 0$ and $2(1+\chi(2c))$ have $x=\pm y$; since $\chi(2) = -1$, exactly one of these two numbers is $4$ and the other is $0$. The remaining $16$ solutions have $x, y \neq 0$ and $x\ne\pm y$, and each pair $1\le A<B\le9$ with $A^{2}+B^{2}\equiv c$ accounts for $8$ of them, namely $(\pm A, \pm B)$ and $(\pm B, \pm A)$. Dividing by $8$ yields two distinct pairs $\{A, B\}$. Hence, there are two pairs for each valid $c$ except $c = 0$ (corresponding to $r = 17$). (You can also verify this statement from Table~\ref{tab:B}, which we obtained by computing $r(A,B)$ for the $36$ pairs.)
\end{proof}

\begin{table}[htbp]
\centering
\caption{The components $B_{r}$ of Theorem~\ref{thm:tenth}, with
coefficients reduced to $\{-9,\dots,9\}$.}\label{tab:B}
\begin{tabular}{rl@{\qquad\qquad}rl}
\toprule
$r$ & $B_{r}$ & $r$ & $B_{r}$\\
\midrule
0 & $G_{3,9}+4yG_{4,6}$ & 10 & $-8G_{1,8}-2G_{4,7}$ \\
1 & $9G_{4,9}+7yG_{2,6}$ & 11 & $9G_{2,7}+2yG_{3,5}$ \\
2 & $-3G_{2,9}+4yG_{6,7}$ & 12 & $-6G_{6,9}+8yG_{4,5}$ \\
3 & $8G_{3,8}-6yG_{7,9}$ & 13 & $3yG_{2,5}+4y^{3}G_{1,3}$ \\
4 & $9G_{4,8}+2yG_{5,6}$ & 14 & $7G_{5,7}+9y^{2}G_{1,4}$ \\
5 & $8G_{2,8}+2G_{5,9}$ & 15 & $-4G_{6,8}-y^{4}G_{1,2}$ \\
6 & $3yG_{1,6}-7yG_{7,8}$ & 16 & $-G_{1,7}+5yG_{8,9}$ \\
7 & $6G_{1,9}+5y^{2}G_{3,4}$ & 17 & $0$ \\
8 & $6G_{5,8}-5y^{3}G_{2,3}$ & 18 & $-5G_{3,6}+yG_{1,5}$ \\
9 & $7G_{3,7}-3y^{2}G_{2,4}$ & & \\
\bottomrule
\end{tabular}
\end{table}

\begin{remark}[The nine products with a repeated factor]\label{rmk:repeated}
$G_{A,B}$ has a repeated factor if and only if $\{A,B\} = \{a, \rho(2a)\}$ for some $a$. There are nine pairs of this form (one for each $a$), which satisfy the cycle of the doubling map $a \mapsto \rho(2a)$. All the other $G_{A,B}$ are products of four distinct $P_{a}$.
\end{remark}

\begin{remark}[Doubling the indices]\label{rmk:doubling-rows}
Replacing $(A,B)$ by $(2A,2B)$ scales $A^{2}+B^{2}$ by $4$, so from \eqref{eq:r-of-pair} it maps the class $r$ to the class $4r+6$. This map fixes $r=17$ only. Since $4$ is a square modulo $19$ of order $9$, it partitions the other eighteen classes into two $9$-cycles, depending upon whether $r+2$ is a quadratic residue or not; these are the two cases used in the proof of Corollary~\ref{cor:two-terms}.
In Section~\ref{sec:exceptional}, we show that the doubling action works also on the coefficients (Remark~\ref{rmk:two-orbits}). Thus, the eighteen nonzero rows of Table~\ref{tab:B} are determined by two of them. The vanishing of $B_{17}$ is not a consequence of this symmetry; it comes from $-1$ being a non-residue modulo $19$.
\end{remark}

\begin{remark}[Comparison with the moduli $11$, $13$ and $17$]
\label{rmk:contrast}
The count in the proof of Corollary~\ref{cor:two-terms} is not specific to $19$. Thanks to \eqref{eq:r-of-pair}, we can compute the analogous count for any prime $p \ge 5$. For each prime $p \ge 5$, we consider pairs $\{A,B\}$ with $1 \le A < B < p/2$, leading to residue classes $r$ where $24r + 10 \equiv 36(A^{2}+B^{2}) \pmod p$. The number of pairs $(A,B)$ such that $A^{2}+B^{2} \equiv c \pmod p$ is given by
\[
\nu_{p}(c) = \frac{p - \chi(-1) - 4 - 2\chi(c)\bigl(1 + \chi(2)\bigr)}{8} \quad (c \not\equiv 0 \pmod p),
\]
and
\[
\nu_{p}(0) = \frac{(p-1)\bigl(1 + \chi(-1)\bigr)}{8},
\]
where $\chi$ is the Legendre symbol on $(\mathbb{Z}/p\mathbb{Z})^{\times}$. This gives
\begin{center}
\begin{tabular}{rl}
\toprule
$p$ & number of quartic products in the components of $f(q)^{10}$\\
\midrule
$11$ & one in each of ten classes, none in one class\\
$13$ & one in each of twelve classes, three in one class\\
$17$ & one in each of eight classes, two in each of eight, four in one\\
$19$ & two in each of eighteen classes, none in one class\\
\bottomrule
\end{tabular}
\end{center}
At $11$, the ten single products are those found by Hirschhorn \cite[(13.4.8)--(13.4.10)]{Hirschhorn2017}, and the empty class corresponds to Winquist's proof of $p(11n+6) \equiv 0$ \cite{Winquist1969}. The modulus $19$ is the first in which no component of $f(q)^{10}$ is given by a single quartic product. We do not know of a reduction of the two-term sums of Table~\ref{tab:B} to single products of the $P_{a}$; none is used below.
\end{remark}

\section{Nine colors}\label{sec:nine}

\begin{lemma}\label{lem:reduction}
Let $p$ be a prime and $1\le k\le p-1$. Then
\[
\sum_{n\ge0}p_{p-k}(n)\,q^{n}\equiv\frac{f(q)^{k}}{f(q^{p})}\pmod p .
\]
\end{lemma}

\begin{proof}
$(1-q^{n})^{p}\equiv1-q^{pn}\pmod p$, so $f(q)^{p}\equiv f(q^{p})$, and
$f(q)^{-(p-k)}=f(q)^{k}/f(q)^{p}$.
\end{proof}

This means that modulo $19$, the colored partition numbers $p_{19-k}(n)$ are governed by the $19$-dissection of $f(q)^{k}$. In this paper, we study the three exponents $k=10$, $k=3$, and $k=2$, corresponding to nine, sixteen, and seventeen colors, respectively.

Let $\Psi_{r}(y)=\sum_{n\ge0}p_{9}(19n+r)\,y^{n}$ for $0\le r\le18$.

\begin{theorem}[All nine-colored components]\label{thm:colored}
For $0\le r\le18$,
\begin{equation}
\Psi_{r} \equiv \frac{P_{0}\,B_{r}}{\mathcal{D}} = \frac{P_{0}\,B_{r}}{P_{1}P_{2}\cdots P_{9}} \pmod{19},
\end{equation}
with $B_{r}$ as in Theorem~\ref{thm:tenth}. Explicitly, $\Psi_{r}/P_{0}$ is congruent to the entry given in Table~\ref{tab:Psi}, which is a sum of two quotients of products of the $P_{a}$ for $r \neq 17$, and zero for $r = 17$.
\end{theorem}

\begin{proof}
From Lemma~\ref{lem:reduction}, we have for $p=19$ and $k=10$ that
\[
\sum_{r} q^{r} \Psi_{r}(y) \equiv \frac{f(q)^{10}}{f(y)} \pmod{19}.
\]
By \eqref{eq:f-of-y} and Theorem~\ref{thm:tenth}, the right side can be rewritten as
\[
\sum_{r} q^{r} \frac{P_{0}^{2}B_{r}}{P_{0}\mathcal{D}}.
\]
Both sides are expressed in the form of \eqref{eq:dissection}, so their respective components must agree. In this expansion, each factor $G_{A,B}/\mathcal{D}$ simplifies as $\mathcal{D}$ divides out the appropriate factors from $G_{A,B}$, leaving exactly one factor in the numerator in the nine cases described in Remark~\ref{rmk:repeated}.
\end{proof}

\begin{corollary}\label{cor:p9}
For every $n\ge0$,
\[
p_{9}(19n+17)\equiv0\pmod{19}.
\]
\end{corollary}

This congruence is known. It is among the results of Kiming and Olsson \cite[Theorem~4]{KimingOlsson1992}, where $\ell=19$ is a special case of \cite[Theorem~2.3(2)]{Boylan2004} and \cite[Theorem~1.2(1)]{LocusWagner2017}, with further consequences discussed here in Section~\ref{sec:family}. What Theorem~\ref{thm:colored} adds is the explicit determination of all the other eighteen components.

\begin{table}[!htbp]
\centering
\caption{The nine-colored components: $\Psi_{r}/P_{0}$ modulo $19$.}
\label{tab:Psi}

\renewcommand{\arraystretch}{2.3}
\begin{tabular}{rl}
\toprule
$r$ & $\Psi_{r}/P_{0}$\\
\midrule
0 & \rule[-11pt]{0pt}{27pt}$\dfrac{1}{P_{1}P_{2}P_{4}P_{5}P_{8}}+\dfrac{4y}{P_{1}P_{3}P_{5}P_{7}P_{8}}$ \\
1 & \rule[-11pt]{0pt}{27pt}$\dfrac{9}{P_{1}P_{2}P_{3}P_{7}P_{8}}+\dfrac{7y}{P_{1}P_{3}P_{5}P_{7}P_{9}}$ \\
2 & \rule[-11pt]{0pt}{27pt}$-\dfrac{3}{P_{1}P_{3}P_{4}P_{5}P_{6}}+\dfrac{4yP_{6}}{P_{2}P_{3}P_{4}P_{5}P_{8}P_{9}}$ \\
3 & \rule[-11pt]{0pt}{27pt}$\dfrac{8P_{8}}{P_{1}P_{2}P_{4}P_{6}P_{7}P_{9}}-\dfrac{6y}{P_{1}P_{4}P_{5}P_{6}P_{8}}$ \\
4 & \rule[-11pt]{0pt}{27pt}$\dfrac{9P_{4}}{P_{1}P_{2}P_{3}P_{5}P_{6}P_{9}}+\dfrac{2y}{P_{2}P_{3}P_{4}P_{7}P_{9}}$ \\
5 & \rule[-11pt]{0pt}{27pt}$\dfrac{8}{P_{1}P_{3}P_{4}P_{5}P_{7}}+\dfrac{2P_{5}}{P_{1}P_{2}P_{3}P_{6}P_{7}P_{8}}$ \\
6 & \rule[-11pt]{0pt}{27pt}$\dfrac{3y}{P_{2}P_{3}P_{4}P_{8}P_{9}}-\dfrac{7y}{P_{2}P_{3}P_{5}P_{6}P_{9}}$ \\
7 & \rule[-11pt]{0pt}{27pt}$\dfrac{6P_{9}}{P_{2}P_{3}P_{4}P_{5}P_{6}P_{7}}+\dfrac{5y^{2}}{P_{2}P_{5}P_{6}P_{8}P_{9}}$ \\
8 & \rule[-11pt]{0pt}{27pt}$\dfrac{6}{P_{1}P_{2}P_{4}P_{7}P_{9}}-\dfrac{5y^{3}}{P_{4}P_{6}P_{7}P_{8}P_{9}}$ \\
9 & \rule[-11pt]{0pt}{27pt}$\dfrac{7}{P_{1}P_{2}P_{5}P_{6}P_{8}}-\dfrac{3y^{2}P_{2}}{P_{1}P_{3}P_{5}P_{7}P_{8}P_{9}}$ \\
10 & \rule[-11pt]{0pt}{27pt}$-\dfrac{8}{P_{2}P_{3}P_{4}P_{5}P_{6}}-\dfrac{2}{P_{1}P_{2}P_{5}P_{6}P_{9}}$ \\
11 & \rule[-11pt]{0pt}{27pt}$\dfrac{9}{P_{1}P_{3}P_{4}P_{6}P_{8}}+\dfrac{2y}{P_{1}P_{4}P_{6}P_{7}P_{9}}$ \\
12 & \rule[-11pt]{0pt}{27pt}$-\dfrac{6}{P_{1}P_{2}P_{5}P_{7}P_{8}}+\dfrac{8y}{P_{2}P_{3}P_{6}P_{7}P_{8}}$ \\
13 & \rule[-11pt]{0pt}{27pt}$\dfrac{3y}{P_{1}P_{4}P_{6}P_{8}P_{9}}+\dfrac{4y^{3}}{P_{5}P_{6}P_{7}P_{8}P_{9}}$ \\
14 & \rule[-11pt]{0pt}{27pt}$\dfrac{7P_{7}}{P_{1}P_{3}P_{4}P_{6}P_{8}P_{9}}+\dfrac{9y^{2}}{P_{2}P_{6}P_{7}P_{8}P_{9}}$ \\
15 & \rule[-11pt]{0pt}{27pt}$-\dfrac{4}{P_{1}P_{3}P_{4}P_{7}P_{9}}-\dfrac{y^{4}P_{1}}{P_{4}P_{5}P_{6}P_{7}P_{8}P_{9}}$ \\
16 & \rule[-11pt]{0pt}{27pt}$-\dfrac{1}{P_{2}P_{3}P_{4}P_{5}P_{9}}+\dfrac{5y}{P_{3}P_{4}P_{5}P_{6}P_{7}}$ \\
17 & \rule[-11pt]{0pt}{27pt}$0$ \\
18 & \rule[-11pt]{0pt}{27pt}$-\dfrac{5P_{3}}{P_{1}P_{2}P_{4}P_{5}P_{7}P_{8}}+\dfrac{y}{P_{2}P_{3}P_{7}P_{8}P_{9}}$ \\
\bottomrule
\end{tabular}
\end{table}

\begin{remark}\label{rmk:first-values}
The constant terms give a quick check on the signs. From
Table~\ref{tab:Psi}, $\Psi_{0}(0)=1$, $\Psi_{1}(0)\equiv9$ and
$\Psi_{2}(0)\equiv-3$, in agreement with $p_{9}(0)=1$, $p_{9}(1)=9$ and
$p_{9}(2)=54$.
\end{remark}

\section{The vanishing class as one member of a family}\label{sec:family}

Corollary~\ref{cor:p9} does not need the dissection. It follows from
\eqref{eq:winquist-f10} alone, for every prime in the same residue class
modulo $4$, and in a sharper form.

\begin{theorem}[Newman; Cooper, Hirschhorn and Lewis]\label{thm:family}
Let $p\ge7$ be a prime with $p\equiv3\pmod4$, and let
$\delta_{p}=5(p^{2}-1)/12$. Write $f(q)^{10}=\sum_{n\ge0}c_{10}(n)q^{n}$.
If $n\equiv\delta_{p}\pmod p$, then
\[
c_{10}(n)=\begin{cases}
p^{4}\,c_{10}(s), & n=p^{2}s+\delta_{p}\ \text{with}\ s\ge0,\\[2pt]
0, & \text{otherwise.}
\end{cases}
\]
Equivalently,
$\sum_{n\equiv\delta_{p}\,(p)}c_{10}(n)q^{n}=p^{4}\,q^{\delta_{p}}
f(q^{p^{2}})^{10}$.
\end{theorem}

\begin{proof}
Let $(a,b)$ index a term of \eqref{eq:winquist-f10} where $N \equiv \delta_{p} \pmod p$. By \eqref{eq:norm-10},
\begin{equation}
X^{2}+Y^{2} = 24N+10 \equiv 24\delta_{p}+10 \equiv 10p^{2} \equiv 0 \pmod p.
\end{equation}
Since $-1$ is a quadratic non-residue modulo $p$ for $p \equiv 3 \pmod 4$, $X$ and $Y$ must both be divisible by $p$.

The multiples of $p$ among the numbers $6a+3$ are of the form $pX'$ with $X' = 6a'+3$, obtained by setting $a = p a' + (p-1)/2$. The sign transforms as $(-1)^{a} = -(-1)^{a'}$. Let $\epsilon = \pm 1$ such that $\epsilon p \equiv 1 \pmod 6$. The multiples of $p$ among the numbers $6b+1$ are of the form $\epsilon p Y'$ with $Y' = 6b'+1$, where $b = \epsilon p b' + (\epsilon p - 1)/6$, so that $(-1)^{b} = (-1)^{b'} (-1)^{(\epsilon p - 1)/6}$. Therefore,
\[
(-1)^{a+b}XY(X^{2}-Y^{2}) = -\epsilon \, (-1)^{(\epsilon p - 1)/6} \, p^{4} \, (-1)^{a'+b'} X' Y' (X'^{2} - Y'^{2}).
\]
When $p \equiv 7 \pmod{12}$, $\epsilon = 1$ and $(p-1)/6$ is odd; when $p \equiv 11 \pmod{12}$, $\epsilon = -1$ and $(p+1)/6$ is even. In both cases, the overall sign factor is $+1$. Finally,
\[
24N+10 = p^{2}(X'^{2}+Y'^{2}) = p^{2}(24N'+10),
\]
where $N'$ is the exponent associated with $(a',b')$, giving $N = p^{2}N' + \delta_{p}$. Thus $(a',b')\mapsto(a,b)$ is a bijection from all terms of \eqref{eq:winquist-f10} onto the terms whose exponent lies in the class of $\delta_{p}$; it multiplies the coefficient by $p^{4}$ and sends $q^{N'}$ to $q^{p^{2}N'+\delta_{p}}$. Dividing both sides by $48$ then yields the theorem.
\end{proof}

\begin{corollary}\label{cor:family}
Let $p\ge11$ be a prime with $p\equiv3\pmod4$. Then
\[
p_{p-10}(N)\equiv0\pmod p\qquad\text{whenever}\quad
N\equiv\frac{5(p^{2}-1)}{12}\pmod p .
\]
For $p=11,19,23,31$ these are
\[
p(11n+6)\equiv0,\qquad p_{9}(19n+17)\equiv0,\qquad
p_{13}(23n+13)\equiv0,\qquad p_{21}(31n+28)\equiv0,
\]
modulo $11$, $19$, $23$ and $31$ respectively.
\end{corollary}

\begin{proof}
By Lemma~\ref{lem:reduction} with $k=10$, the generating function of $p_{p-10}$ satisfies
\begin{equation}
\sum_{n\ge0}p_{p-10}(n)q^{n} \equiv \frac{f(q)^{10}}{f(q^{p})} \pmod p.
\end{equation}
Multiplication preserves residue classes of exponents, and the latter factor is a series in $q^{p}$. Theorem~\ref{thm:family} states that the part of $f(q)^{10}$ in the class of $\delta_{p}$ is divisible by $p^{4}$.
\end{proof}

\begin{remark}[What is known]\label{rmk:family-known}
Nothing in Theorem~\ref{thm:family} or Corollary~\ref{cor:family} is new. For $p\equiv11\pmod{12}$, Theorem~\ref{thm:family} is a case of Newman's theorem \cite{Newman1955}, proved using modular functions. Newman stated the case $p\equiv7\pmod{12}$, which contains $p=19$, without proof \cite[(7)]{Newman1957a}. Cooper, Hirschhorn, and Lewis proved both cases \cite[Theorem~1]{CHL2000} in the form
\begin{equation}
c_{10}\left(pn+\frac{5(p^{2}-1)}{12}\right) = p^{4}c_{10}\left(\frac{n}{p}\right),
\end{equation}
with $c_{10}(n/p)=0$ when $p\nmid n$; putting $n=ps$ gives the statement above. Their proof is the argument given here: both indices of Winquist's series are divisible by $p$, and one rescales and follows the signs. 

Newman deduced congruences for colored partitions from such exact relations as in the proof of Corollary~\ref{cor:family} \cite[Theorem~1]{Newman1957a}; his hypotheses give $p\equiv11\pmod{12}$ when there are ten factors. The case $p=11$ of the corollary is Ramanujan's congruence, and the argument is Winquist's \cite{Winquist1969}. The function $\eta(12\tau)^{10} = q^{5}f(q^{12})^{10}$ is a cusp form of weight $5$ with complex multiplication by $\mathbb{Q}(i)$, which is why it is lacunary \cite{Serre1985}; a prime $p\equiv3\pmod4$ is inert in $\mathbb{Q}(i)$, and the relation $c_{10}(p^{2}s+\delta_{p})=p^{4}c_{10}(s)$ is the Hecke relation at such a prime, with the factor $p^{4}=p^{k-1}$ coming from the weight. 
Corollary~\ref{cor:family} is \cite[Theorem~2.3(2)]{Boylan2004} with $i=10$, which Boylan describes as well known to experts and proves from Serre's decomposition \cite[\S2.5]{Serre1985} of $\eta(12\tau)^{10}$ into two CM newforms; it is also \cite[Theorem~1.2(1)]{LocusWagner2017} with $r=10$ and $s=1$. At $p=19$, it appears in the list of Kiming and Olsson \cite[Theorem~4]{KimingOlsson1992}.
The above proof is included in the paper since it is brief and is written in the notation of Section~\ref{sec:prelim}. The power $p^{4}$ does not pass to $p_{p-10}$ because Lemma~\ref{lem:reduction} holds modulo $p$ only.
\end{remark}

\section{Seventeen and sixteen colors}\label{sec:seventeen}

\subsection{The square of Euler's product}

By Euler's pentagonal number theorem,
\begin{equation}\label{eq:f-squared}
f(q)^{2}=\sum_{m,n\in\mathbb{Z}}(-1)^{m+n}
q^{\frac{m(3m-1)}{2}+\frac{n(3n-1)}{2}},
\end{equation}
and with $X=6m-1$, $Y=6n-1$ the exponent $N$ satisfies
\begin{equation}\label{eq:norm-2}
24N+2=X^{2}+Y^{2}.
\end{equation}
More generally, for a polynomial $H(X,Y)$ put
\[
S_{H}(q)=\sum_{m,n\in\mathbb{Z}}(-1)^{m+n}H(6m-1,6n-1)\,
q^{\frac{m(3m-1)}{2}+\frac{n(3n-1)}{2}},
\]
so that $S_{1}=f(q)^{2}$.

\begin{theorem}\label{thm:square}
Let $p\ge7$ be a prime with $p\equiv3\pmod4$, let
$\delta=(p^{2}-1)/12$, and let $H(X,Y)$ be homogeneous of degree $h$ and
even in each variable. Write $S_{H}=\sum_{n\ge0}c_{H}(n)q^{n}$. If
$n\equiv\delta\pmod p$, then
\[
c_{H}(n)=\begin{cases}
p^{h}\,c_{H}(s), & n=p^{2}s+\delta\ \text{with}\ s\ge0,\\[2pt]
0, & \text{otherwise.}
\end{cases}
\]
In particular
$\sum_{n\equiv\delta\,(p)}[q^{n}]f(q)^{2}\;q^{n}=q^{\delta}f(q^{p^{2}})^{2}$.
\end{theorem}

\begin{proof}
If $N \equiv \delta \pmod p$, by \eqref{eq:norm-2} we have $X^{2}+Y^{2} \equiv 24\delta+2 \equiv 2p^{2} \equiv 0 \pmod p$, and so $p$ divides both $X$ and $Y$. As before, letting $\epsilon p \equiv 1 \pmod 6$, the multiples of $p$ in the quantity $6m-1$ are given by $\epsilon p X'$, where $X'$ is set to $6m'-1$; consequently, $m = \epsilon p m' - (\epsilon p - 1)/6$, and $(-1)^{m} = (-1)^{m'} (-1)^{(\epsilon p - 1)/6}$. The same holds for $n$. 
Since the additional sign is the same for $m$ and for $n$, we have the sign relation $(-1)^{m+n} = (-1)^{m'+n'}$. Since $H$ is even in each variable, $H(\epsilon p X', \epsilon p Y') = p^{h} H(X', Y')$. Finally, the relation $24N+2 = p^{2}(24N'+2)$ yields $N = p^{2}N' + \delta$.
\end{proof}

\begin{remark}\label{rmk:square-known}
As with Theorem~\ref{thm:family}, Theorem~\ref{thm:square} is classical in substance. For $H=1$, it is the case of two factors of \cite[Theorem~1]{CHL2000}, proved there for $p \equiv 5, 7,$ or $11 \pmod{12}$ by the argument used above; it goes back to Newman, who proved it for $p \equiv 11 \pmod{12}$ \cite[Theorem~2]{Newman1953} and stated it for $p \not\equiv 1 \pmod{12}$ \cite[(4)]{Newman1957a}. For $H=H_{4}$ and $H=H_{8}$, the same argument holds true when the harmonic weight is included. 
If $H$ is the real part of $(X+iY)^{h}$, then the series
\begin{equation}
q\,S_{H}(q^{12}) = \sum_{m,n\in\mathbb{Z}} (-1)^{m+n} H(X,Y) \, q^{(X^{2}+Y^{2})/2}, \qquad X = 6m-1, \ Y = 6n-1,
\end{equation}
is a theta series of $\mathbb{Q}(i)$ with a Hecke character, representing a cusp form of weight $h+1$, and the theorem is its Hecke relation at primes inert in $\mathbb{Q}(i)$. For $H=1$, this reduces to the lacunarity of $\eta(12\tau)^{2}$ \cite{Serre1985}. Huber, Liu, McLaughlin, Ye, Yuan, and Zhang \cite[Theorem~1.1(14)]{HuberEtAl2023} characterize exactly when its coefficients vanish: $\eta(12\tau)^{2}$ has complex multiplication by both $\mathbb{Q}(i)$ and $\mathbb{Q}(\sqrt{-3})$, and the coefficient of $q^{n}$ is zero whenever a prime $p \not\equiv 1 \pmod{12}$ divides $n$ to an odd power. 
The vanishing part of Theorem~\ref{thm:square} for $H=1$ therefore holds for all primes $p \not\equiv 1 \pmod{12}$; the elementary argument above covers those with $p \equiv 3 \pmod 4$ and gives the exact value on the sublattice. For $H=H_{4}$, Proposition~\ref{prop:eisenstein} identifies $E_{4}(12\tau)\eta(12\tau)^{2}$ with a theta series of this kind; in Serre's decomposition \cite[\S2.5]{Serre1985}, recalled in \cite[\S5]{Boylan2004}, it is half the sum of the two CM newforms whose difference is $96\,\eta(12\tau)^{10}$.
\end{remark}

\begin{corollary}[Seventeen colors]\label{cor:seventeen}
\[
\sum_{n\ge0}p_{17}(19n+11)\,y^{n}\equiv\frac{y\,f(y^{19})^{2}}{f(y)}
\pmod{19}.
\]
In particular $p_{17}(11)\equiv0$ and $p_{17}(19n+11)\equiv p(n-1)$ for
$1\le n\le19$. More generally, if $p\ge7$ is a prime with $p\equiv3\pmod4$
and $(p^{2}-1)/12=p\,\kappa_{p}+r_{p}$ with $0\le r_{p}<p$, then
\[
\sum_{n\ge0}p_{p-2}(pn+r_{p})\,x^{n}\equiv
\frac{x^{\kappa_{p}}f(x^{p})^{2}}{f(x)}\pmod p .
\]
\end{corollary}

\begin{proof}
From Lemma~\ref{lem:reduction} with $k=2$, the generating function of $p_{p-2}$ is congruent to $f(q)^{2}/f(q^{p})$ modulo $p$. Taking the part with exponents in the class of $\delta$ and using Theorem~\ref{thm:square} with $H=1$ gives
\begin{equation}
\frac{q^{\delta}f(q^{p^2})^{2}}{f(q^{p})},
\end{equation}
and then one sets $x = q^{p}$. For $p=19$, $\delta = 30 = 19+11$. That is, the assertion for $p_{17}(19n+11) \equiv p(n-1) \pmod{19}$ holds true since $f(y^{19})^{2} = 1 + O(y^{19})$.
\end{proof}

In the case $p \equiv 5 \pmod{12}$, we have from the case $r=2$ of \cite[Theorem~1]{CHL2000} the equality with the sign $-1$ (and an identical deduction gives)
\begin{equation}
\sum_{n\ge0} p_{p-2}(pn+r_{p})\,x^{n} \equiv -\frac{x^{\kappa_{p}}f(x^{p})^{2}}{f(x)} \pmod p,
\end{equation}
so, up to the sign, we have our formula for every prime $p \ge 5$ with $p \not\equiv 1 \pmod{12}$. The proof of Theorem~\ref{thm:square} does not reach primes with $p \equiv 1 \pmod 4$ because in that case the congruence $X^{2}+Y^{2} \equiv 0 \pmod p$ does not imply that $p$ divides $X$.

\subsection{Two Eisenstein weights}

Let $P$, $Q$, $R$ be Ramanujan's series
\[
P=1-24\sum_{n\ge1}\sigma_{1}(n)q^{n},\qquad
Q=1+240\sum_{n\ge1}\sigma_{3}(n)q^{n},\qquad
R=1-504\sum_{n\ge1}\sigma_{5}(n)q^{n},
\]
and let $\vartheta=24\,q\,\frac{d}{dq}$. Ramanujan's differential equations
\cite{Ramanujan1916} read
\begin{equation}\label{eq:ramanujan-ode}
\vartheta P=2(P^{2}-Q),\qquad\vartheta Q=8(PQ-R),\qquad
\vartheta R=12(PR-Q^{2}),
\end{equation}
and $Q^{2}=1+480\sum_{n\ge1}\sigma_{7}(n)q^{n}$.

\begin{lemma}\label{lem:U}
For $k\ge0$ let $U_{k}=f(q)^{-1}\sum_{m\in\mathbb{Z}}(-1)^{m}(6m-1)^{2k}
q^{m(3m-1)/2}$. Then $U_{0}=1$, $U_{k+1}=\vartheta U_{k}+P\,U_{k}$, and
\begin{gather*}
U_{1}=P,\qquad U_{2}=3P^{2}-2Q,\qquad U_{3}=15P^{3}-30PQ+16R,\\
U_{4}=105P^{4}-420P^{2}Q+448PR-132Q^{2}.
\end{gather*}
\end{lemma}

\begin{proof}
Since $(\vartheta+1)q^{m(3m-1)/2} = (6m-1)^{2}q^{m(3m-1)/2}$, we have $fU_{k} = (\vartheta+1)^{k}f$. Logarithmic differentiation of the product gives $\vartheta f = -24f\sum\sigma_{1}(n)q^{n}$, that is, $(\vartheta+1)f = Pf$. Hence,
\[
fU_{k+1} = (\vartheta+1)(fU_{k}) = f\,\vartheta U_{k} + U_{k}(\vartheta+1)f = f(\vartheta U_{k} + PU_{k}).
\]
Applying Ramanujan's differential equations \eqref{eq:ramanujan-ode} repeatedly (for instance, $U_{2} = 2(P^{2}-Q) + P^{2}$) yields the four values.
\end{proof}

\begin{proposition}\label{prop:eisenstein}
Let $H_{4}=X^{4}-6X^{2}Y^{2}+Y^{4}$ and
$H_{8}=X^{8}-28X^{6}Y^{2}+70X^{4}Y^{4}-28X^{2}Y^{6}+Y^{8}$, the real parts of
$(X+iY)^{4}$ and $(X+iY)^{8}$. Then
\[
S_{H_{4}}=-4\,f(q)^{2}Q,\qquad S_{H_{8}}=16\,f(q)^{2}Q^{2}.
\]
\end{proposition}

\begin{proof}
For a monomial, $S_{X^{2i}Y^{2j}}=f(q)^{2}U_{i}U_{j}$. Hence
$S_{H_{4}}=f^{2}(2U_{2}-6U_{1}^{2})=f^{2}(6P^{2}-4Q-6P^{2})$, and
$S_{H_{8}}=f^{2}(2U_{4}-56U_{3}U_{1}+70U_{2}^{2})$. In the last bracket the
coefficients of $P^{4}$, $P^{2}Q$ and $PR$ are $210-840+630$,
$-840+1680-840$ and $896-896$, all zero, and the coefficient of $Q^{2}$ is
$-264+280=16$.
\end{proof}

\begin{corollary}[Weighted congruences]\label{cor:weighted}
Let $p\ge7$ be a prime with $p\equiv3\pmod4$ and let
$\delta=(p^{2}-1)/12$. Then the part of $f(q)^{2}Q$ with exponents in the
class of $\delta$ modulo $p$ is $p^{4}q^{\delta}\,(f^{2}Q)(q^{p^{2}})$, the
corresponding part of $f(q)^{2}Q^{2}$ is $p^{8}q^{\delta}\,(f^{2}Q^{2})
(q^{p^{2}})$, and for every $N\equiv\delta\pmod p$
\begin{align*}
p_{p-2}(N)+240\sum_{m=1}^{N}\sigma_{3}(m)\,p_{p-2}(N-m)&\equiv0\pmod p,\\
p_{p-2}(N)+480\sum_{m=1}^{N}\sigma_{7}(m)\,p_{p-2}(N-m)&\equiv0\pmod p.
\end{align*}
For $p=19$ these read, for $N=19n+11$,
\begin{align*}
p_{17}(N)+12\sum_{m=1}^{N}\sigma_{3}(m)\,p_{17}(N-m)&\equiv0\pmod{19},\\
p_{17}(N)+5\sum_{m=1}^{N}\sigma_{7}(m)\,p_{17}(N-m)&\equiv0\pmod{19}.
\end{align*}
\end{corollary}

\begin{proof}
The first two statements follow from Theorem~\ref{thm:square} for $H_{4}$ and $H_{8}$ combined with Proposition~\ref{prop:eisenstein}: the factors $-4$ and $16$ are prime to $p$. By Lemma~\ref{lem:reduction},
\begin{equation}
Q \sum_{n\ge0} p_{p-2}(n)q^{n} \equiv \frac{f(q)^{2}Q}{f(q^{p})} \pmod p,
\end{equation}
and the part of the right-hand side in the class of $\delta$ is divisible by $p^{4}$. The coefficient of $q^{N}$ on the left is the first displayed sum. The second is obtained in the same way from $Q^{2}$.
\end{proof}

By Corollary~\ref{cor:seventeen}, a whole component is determined rather than just a vanishing class. In \cite[Theorem~1]{KimingOlsson1992}, it is shown that the only class that can carry a congruence for $p_{p-2}$ is the one with $24a \equiv p-2 \pmod p$, that is, $a \equiv -1/12 \equiv \delta$, whereas there is no mention of $i=2$ being exceptional in \cite[Theorem~2.3]{Boylan2004}. In fact, $19$ is not exceptional for $r=17$ \cite[Theorem~1.3]{Boylan2004}. Corollary~\ref{cor:seventeen} shows instead what happens: the class is non-zero and appears in closed form.
Closed forms of this kind for particular components of colored partition functions go back to Newman \cite[Theorem~1]{Newman1957b}, who obtained, for instance, $\sum_{n} p(13n+6)x^{n} \equiv 11(x;x)_{\infty}^{11} \pmod{13}$, in the classes $r(p-1)/24$. Corollary~\ref{cor:seventeen} deals with the class $2(p^{2}-1)/24$ instead.
It follows at once from the case $H=1$ of Theorem~\ref{thm:square}, that is, from \cite[Theorem~1]{CHL2000}. However, we have not actually been able to find this exact formula in the existing literature.

\begin{remark}\label{rmk:weighted-scope}
Under the exchange $X \leftrightarrow Y$, the harmonic polynomial $\operatorname{Re}(X+iY)^{h}$ changes sign when $h \equiv 2 \pmod 4$. Thus, $S_{H}$ vanishes identically for such $h$, and the method gives nothing in weights $6$ and $10$.
\end{remark}

\subsection{Sixteen colors}

For completeness we record what Proposition~\ref{prop:jacobi} gives.
Write $t_{a}=19k_{a}+r_{a}$ with $0\le r_{a}\le18$, so that
$(k_{a},r_{a})=(1,17),(1,9),(1,2),(0,15),(0,10),(0,6),(0,3),(0,1),(0,0)$
for $a=1,\dots,9$.

\begin{proposition}\label{prop:sixteen}
For $1\le a\le9$,
\[
\sum_{n\ge0}p_{16}(19n+r_{a})\,y^{n}\equiv
\frac{(-1)^{a}\,2a\;y^{k_{a}}}{\prod_{b\ne a}P_{b}}\pmod{19},
\]
the product being over $1\le b\le9$, $b\ne a$, and
\[
p_{16}(19n+r)\equiv0\pmod{19}\qquad\text{for}\quad
r\in\{4,5,7,8,11,12,13,14,16,18\}.
\]
\end{proposition}

\begin{proof}
By Lemma~\ref{lem:reduction} with $k=3$, Proposition~\ref{prop:jacobi} and
\eqref{eq:f-of-y}, the generating function of $p_{16}$ is congruent to
$\sum_{a}(-1)^{a}2a\,q^{t_{a}}P_{a}/\mathcal{D}$.
\end{proof}

The ten classes are those $r$ for which $8r+1$ is either divisible by $19$
or a non-residue, in agreement with \cite[Theorem~2(2)]{KimingOlsson1992}.

\subsection{Thirteen colors}

Kiming and Olsson's list at $19$ \cite[Theorem~4]{KimingOlsson1992}
consists of $p_{7}(19n+9)$, $p_{9}(19n+17)$ and $p_{13}(19n+14)$. The second
is Corollary~\ref{cor:p9}. The third also falls to
Proposition~\ref{prop:jacobi}.

\begin{proposition}\label{prop:thirteen}
$p_{13}(19n+14)\equiv0\pmod{19}$ for every $n\ge0$.
\end{proposition}

\begin{proof}
Thus, combining Lemma~\ref{lem:reduction} with $k=6$ and Proposition~\ref{prop:jacobi}, we see that the generating function of $p_{13}$ modulo $19$ is congruent to
\begin{equation}
\frac{P_{0}\left(\sum_{a}(-1)^{a}2a\,q^{t_{a}}P_{a}\right)^{2}}{\mathcal{D}},
\end{equation}
whose exponents are of the form $t_{a}+t_{b}+19k$. The residues of $t_{1},\dots,t_{9}$ are $17, 9, 2, 15, 10, 6, 3, 1, 0$, and no two of them, equal or not, add up to $14$ modulo $19$.
\end{proof}

The first congruence, which Boylan calls superexceptional, is not explained by complex multiplication \cite[\S2]{Boylan2004}. We return to it in Section~\ref{sec:questions}.

\section{The Winquist--Jacobi bridge}\label{sec:bridge}

\subsection{Homogeneous variables}

As such, we can extend the definition of $P_{x}$ to all integers $x$ by
\[
P_{x} = (y^{x};y^{19})_{\infty}(y^{19-x};y^{19})_{\infty}.
\]
Then $P_{19-x} = P_{x}$, $P_{x} = 0$ when $x \equiv 0 \pmod{19}$, and as in the proof of Lemma~\ref{lem:winquist-special}, $P_{-x} = -y^{-x}P_{x}$. Extending $t_{x}$ by \eqref{eq:t-and-d} and putting
\[
Z_{x} = q^{t_{x}}P_{x} \qquad (x \in \mathbb{Z}),
\]
we have by \eqref{eq:t-d-shift} and $t_{19-x} = t_{x}$ that
\begin{equation}\label{eq:Z-rules}
Z_{19-x} = Z_{x}, \qquad Z_{-x} = -Z_{x}, \qquad Z_{x+19} = -Z_{x},
\end{equation}
hence every $Z_{x}$ is either $0$, or $\pm Z_{j}$ for some $1 \le j \le 9$. Let $\mathcal{R} = \mathbb{F}_{19}[Z_{1},\dots,Z_{9}]$ with $Z_{1},\dots,Z_{9}$ as indeterminates, and let
\[
\operatorname{ev} \colon \mathcal{R} \to \mathbb{F}_{19}((q))
\]
be the reduced ring homomorphism sending $Z_{j}$ to $q^{t_{j}}P_{j}$ modulo $19$. For a sequence of integers $e_{1},\dots,e_{9}$, we define the monomial $Z_{1}^{e_{1}}\cdots Z_{9}^{e_{9}}$ to have weight $\sum e_{j}t_{j}$ modulo $19$. The evaluation of a monomial of weight $w$ is $q^{w}$ times a power series in $y$, so the evaluation of a polynomial whose monomials all have weight $w$ has its exponents in the class of $w$ modulo $19$.

In Proposition~\ref{prop:jacobi}, $t_{a}$ appears as the exponent of $q$ in front of $P_{a}$, so this proposition reads as $f(q)^{3} \equiv P_{0}\operatorname{ev}(C)$ with
\[
C = \sum_{a=1}^{9} (-1)^{a} \, 2a \, Z_{a}.
\]
For the tenth power, $t_{B+A} + t_{B-A} = 2t_{B}+A^{2}$ gives
\begin{align*}
t_{A} + t_{B} + t_{A+B} + t_{B-A} &= t_{A} + 3t_{B} + A^{2} \\
&= \frac{3A^{2}+3B^{2}-19A-57B+360}{2} = \varepsilon(A,B) + 30.
\end{align*}
Thus, $G_{A,B}$ for each term $q^{\varepsilon(A,B)}$ in Theorem~\ref{thm:tenth} is evaluated as $q^{-30}Z_{A}Z_{B}Z_{A+B}Z_{B-A}$, and $f(q)^{10} \equiv q^{-30}P_{0}^{2}\operatorname{ev}(B)$ with
\[
B = -3 \sum_{1 \le A < B' \le 9} (-1)^{A+B'}AB'\,(B'^{2}-A^{2})\,Z_{A}Z_{B'}Z_{A+B'}Z_{B'-A}.
\]
(We denote by $B'$ the summation variable to avoid confusion with the polynomial $B$.)

\subsection{The bridge}

Let $\Phi_{s}(y) = \sum_{n\ge0} p(19n+s) \, y^{n}$ for $0 \le s \le 18$, and let $H_{0},\dots,H_{18}$ be the components of $\operatorname{ev}(C)^{3}$ in the sense of \eqref{eq:dissection}.

\begin{proposition}[Winquist--Jacobi bridge]\label{thm:bridge}
For every $0 \le t \le 18$,
\begin{equation}
\sum_{s=0}^{18} y^{[s>t]} \, B_{(t-s)\bmod 19} \, \Phi_{s} \equiv P_{0}H_{t} \pmod{19},
\end{equation}
where $[s>t]$ is $1$ if $s>t$ and $0$ otherwise, and $B_{r}$ are the components of Theorem~\ref{thm:tenth}.
\end{proposition}

\begin{proof}
Consider $f(q)^{10} \cdot f(q)^{-1} = f(q)^{9} = \bigl(f(q)^{3}\bigr)^{3}$. Modulo $19$, the left side is $P_{0}^{2}\sum_{r}q^{r}B_{r} \cdot \sum_{s}q^{s}\Phi_{s}$, and the right side is $P_{0}^{3}\sum_{t}q^{t}H_{t}$. Collecting the left side by the class of $r+s$, noting that $q^{r+s} = q^{t}y$ when $r+s = t+19$, and dividing by $P_{0}^{2}$ yields the result.
\end{proof}

\begin{proposition}[The scalar identity]\label{thm:scalar}
Let $F = B \, C^{3} \in \mathcal{R}$, and for each class $w$ let $F_{w}$ be the sum of the monomials of $F$ of weight $w$. Then
\begin{equation}
\operatorname{ev}(F_{11}) \equiv q^{30} \, \frac{P_{1}P_{2}\cdots P_{9}}{P_{0}^{4}}, \qquad \operatorname{ev}(F_{w}) \equiv 0 \quad (w \not\equiv 11 \pmod{19}).
\end{equation}
\end{proposition}

\begin{proof}
As in the proof of Lemma~\ref{lem:reduction} and by \eqref{eq:f-of-y}, we have $f(q)^{19} \equiv f(y) = P_{0}\mathcal{D}$, while $f(q)^{19} = f(q)^{10}(f(q)^{3})^{3} \equiv q^{-30}P_{0}^{5}\operatorname{ev}(F)$. Thus $\operatorname{ev}(F) \equiv q^{30}\mathcal{D}/P_{0}^{4}$, all of whose exponents lie in the class of $30 \equiv 11$. The part of $\operatorname{ev}(F)$ in a class $w$ is $\operatorname{ev}(F_{w})$.
\end{proof}

\begin{remark}\label{rmk:bridge-limits}
Proposition~\ref{thm:scalar} gives, without any computation, eighteen polynomials $F_{w}$ of degree $7$ that vanish on the theta series, and one, $F_{11}$, whose value is a quotient of products. Proposition~\ref{thm:bridge} provides a linear system for $\Phi_{0},\dots,\Phi_{18}$ with known coefficients; eliminating the $B_{r}$ by Proposition~\ref{thm:scalar} returns $f(q)^{-1} = f(q)^{18}/f(y)$, so the system by itself produces no short formula for any individual $\Phi_{s}$. The exceptional component is the one case where the structure of $F_{11}$ can be brought to bear, and this is the subject of the next section.
\end{remark}
\section{The exceptional component}\label{sec:exceptional}

The residue $s=4$ is exceptional because $24\cdot4 \equiv 1 \pmod{19}$. Two earlier facts frame this section. Kolberg \cite[proof of Lemma~4]{Kolberg1959}, building bases of modular functions on $\Gamma_{0}(19)$ out of Fine's functions $W_{k}$ \cite{Fine1956}, uses
\[
f_{2} = W_{1}W_{7}W_{8} + W_{2}W_{5}W_{3} + W_{4}W_{9}W_{6},
\]
whose three index blocks $\{1,7,8\}$, $\{2,3,5\}$, and $\{4,6,9\}$ are the cosets of $\{1,7,8\} = \{\pm 1, \pm 7, \pm 8\}$ in $(\mathbb{Z}/19)^{\times}/\{\pm 1\}$. These blocks are precisely the same blocks that form the base of the cubic $S$ shown below. O'Brien gives an exact formula for this exceptional component over $\mathbb{Z}$ in Theorem~6.1 of his thesis \cite{OBrien1965}.

Building on this exceptional component over $\mathbb{Z}$, O'Brien deduces \cite[(6.59)]{OBrien1965} that
\[
\Phi_{4} \equiv 5f(y^{19})f(y)^{16}S \pmod{19},
\]
with
\[
S = P_{4}P_{6}P_{9} + yP_{1}P_{7}P_{8} - y^{2}P_{2}P_{3}P_{5}.
\]
Since $f(y)^{19} \equiv f(y^{19}) = P_{0}$ and $f(y) = P_{0}\mathcal{D}$, this becomes
\begin{equation}\label{eq:obrien}
\Phi_{4} \equiv \frac{5\,S}{P_{0}\,\mathcal{D}^{3}} \pmod{19}.
\end{equation}

His argument contains more, and we summarize what is used below. He starts from the $19$-dissection of $f(q)$ itself in terms of Fine's functions, whose indices run along the doubling cycle, and works throughout with the symmetry that interchanges his variables cyclically (index doubling). The products he forms over the cosets of $\{1,7,8\}$ are what lead to the blocks. Setting
\[
m_{1} = yP_{1}P_{7}P_{8}, \qquad m_{2} = -y^{2}P_{2}P_{3}P_{5}, \qquad m_{3} = P_{4}P_{6}P_{9},
\]
so that $S = m_{1}+m_{2}+m_{3}$, he proves \cite[(6.47), (6.51), (6.52)]{OBrien1965} exactly that
\begin{gather*}
m_{1}m_{2}m_{3} = -y^{3}\mathcal{D}, \qquad
m_{1}m_{2} + m_{2}m_{3} + m_{3}m_{1} = y\mathcal{D}^{2}, \\
m_{1}m_{2}^{2} + m_{2}m_{3}^{2} + m_{3}m_{1}^{2} = -2y^{3}\mathcal{D},
\end{gather*}
and then, as part of \cite[(6.57) and Theorem~6.1]{OBrien1965}, establishes the relation which in the present notation reads
\begin{equation}\label{eq:obrien-cubic-Z}
S^{3} = 8y\mathcal{D}^{2}S + \mathcal{D}^{5} + 19y^{3}\mathcal{D},
\end{equation}
an identity over $\mathbb{Z}$. (In his notation, $\varepsilon = yS/\mathcal{D}^{3}$, $F = y^{-3}\mathcal{D}^{4}$, and the relation is $\varepsilon^{3} = (8\varepsilon+1)/F + 19/F^{2}$.) Modulo $19$, this yields
\begin{equation}\label{eq:obrien-cubic}
S^{3} \equiv 8y\mathcal{D}^{2}S + \mathcal{D}^{5} \pmod{19}.
\end{equation}

Kolberg had earlier given four congruences for the same component \cite[(3), (6)--(8)]{Kolberg1959}; in his notation, where $\varphi(x)^{r} = \sum p_{r}(n)x^{n}$ for every integer $r$, the first is $p(19n+4) \equiv p_{71}(19n+1) \pmod{19}$.

\begin{remark}[Cubic is a model for $X_{0}(19)$]\label{rmk:weierstrass}
Basis functions for Kolberg at $q=19$ \cite[proof of Lemma~4]{Kolberg1959} are, in the present notation,
\[f_{2} = \frac{S\,\mathcal{D}}{y^{2}} + 5, \qquad f_{3} = w = \frac{f(y)^{4}}{y^{3}f(y^{19})^{4}} = \frac{\mathcal{D}^{4}}{y^{3}};\]
the second is \eqref{eq:f-of-y}. For the first, the definition of $W_{k}$ and the rule $P_{x+19} = -y^{-x}P_{x}$ show that the three blocks of $f_{2}$ are $-m_{3}/m_{2}$, $-m_{1}/m_{3}$ and $-m_{2}/m_{1}$, so that $f_{2} = -(m_{1}^{2}m_{2} + m_{2}^{2}m_{3} + m_{3}^{2}m_{1})/(m_{1}m_{2}m_{3})$, and by O'Brien's three identities,
\[m_{1}^{2}m_{2} + m_{2}^{2}m_{3} + m_{3}^{2}m_{1} = S \cdot y\mathcal{D}^{2} - 3(-y^{3}\mathcal{D}) - (-2y^{3}\mathcal{D}) = yS\mathcal{D}^{2} + 5y^{3}\mathcal{D},\]
using $\sum_{i \ne j} m_{i}^{2}m_{j} = e_{1}e_{2} - 3e_{3}$ for the elementary symmetric functions of $m_{1}, m_{2}, m_{3}$ where $e_{1}=S$. 

Putting $X = f_{2}-5$ and $Y = f_{3}$, multiplying \eqref{eq:obrien-cubic-Z} by $\mathcal{D}^{3}/y^{6}$ gives
\[Y^{2} + 8XY + 19Y = X^{3}.\]
Kolberg shows \cite{Kolberg1959} that $f_{2}$ and $f_{3}$ are modular functions on $\Gamma_{0}(19)$ with integral coefficients, holomorphic away from the cusp at infinity, where they have poles of orders $2$ and $3$. Since $X_{0}(19)$ has genus $1$, this means that the cubic is a Weierstrass model for some elliptic curve $E$ with discriminant $-19^{3}$ and $c_{4}=448$, yielding $X_{0}(19)$. The linear change of coordinates $X = x-5$, $Y = y-4x+11$ transforms it to the minimal model of $X_{0}(19)$, Cremona's curve $19a1$ \cite{Cremona1997}:
\[y^{2} + y = x^{3} + x^{2} - 9x - 15.\]
In this model, $x = f_{2}$ is Kolberg's function itself, and $y = f_{3} + 4f_{2} - 11$. Thus, O'Brien's cubic \eqref{eq:obrien-cubic-Z} is expressed in Kolberg's coordinates. Modulo $19$, this gives $Y^{2} + 8XY = X^{3}$, which has a node at the origin since $19 \nmid c_{4}$. Below we provide that equation rewritten as \eqref{eq:obrien-factor}.
\end{remark}

The cubic $S$, its blocks, the doubling symmetry behind them, the formula \eqref{eq:obrien} and the relation \eqref{eq:obrien-cubic-Z} are all due to O'Brien. By \eqref{eq:obrien-cubic},
\begin{equation}\label{eq:obrien-factor}
S \cdot \frac{S^{2}-8y\mathcal{D}^{2}}{f(y)^{4}} \equiv \frac{\mathcal{D}}{P_{0}^{4}} \pmod{19}.
\end{equation}
Since $S(0)=1$, $S$ is a unit in $\mathbb{F}_{19}[[y]]$. As a statement of divisibility, \eqref{eq:obrien-factor} says nothing new; its true content lies in the explicit cofactor. The right-hand side is, up to a power of $q$, the value of the weight-$11$ part via Proposition~\ref{thm:scalar}. 

We show below that a doubling automorphism of the polynomial ring of the theta functions, which is O'Brien's symmetry in polynomial form, preserves the three-term relations, multiplies the weight-$11$ part by $4$, and fixes the cubic. This is compatible with a factorization in the polynomial ring but does not imply one.
\subsection{The three-term relations}

For integers $a, b, c, d$, put
\[
\Theta(a,b,c,d) = Z_{a+d}Z_{a-d}Z_{b+c}Z_{b-c} - Z_{a+c}Z_{a-c}Z_{b+d}Z_{b-d} + Z_{a+b}Z_{a-b}Z_{c+d}Z_{c-d},
\]
where $Z$ is evaluated in $\mathcal{R}$ via the composition rules of \eqref{eq:Z-rules}. Let $\mathcal{I}_{\Theta} \subset \mathcal{R}$ be the ideal generated by these elements.

\begin{lemma}\label{lem:theta}
$\operatorname{ev}(\Theta(a,b,c,d)) = 0$ for all $a,b,c,d$, so $\operatorname{ev}$ vanishes on $\mathcal{I}_{\Theta}$. Each of the three monomials of $\Theta(a,b,c,d)$ has weight congruent to $9 + a^{2} + b^{2} + c^{2} + d^{2} \pmod{19}$.
\end{lemma}

\begin{proof}
We employ Weierstrass's three-term identity for the theta function
\[
\vartheta(z) = (z; y^{19})_{\infty}(y^{19}/z; y^{19})_{\infty},
\]
in multiplicative form \cite[(2.7)]{Koornwinder2014}. This is equivalent to Jacobi's five-term identity \cite{Koornwinder2014}. Writing $\vartheta(y^{x}) = P_{x}$ and scaling by the appropriate power of $y$, the identity becomes
\[
P_{a+d}P_{a-d}P_{b+c}P_{b-c} - P_{a+c}P_{a-c}P_{b+d}P_{b-d} + y^{b-c} P_{a+b}P_{a-b}P_{c+d}P_{c-d} = 0.
\]
The exponent relation $t_{u+v} + t_{u-v} = 2t_{u} + v^{2}$ shows that the first two terms have $q$-weight $2t_{a} + 2t_{b} + c^{2} + d^{2}$, while the third term has $q$-weight $2t_{a} + 2t_{c} + b^{2} + d^{2}$. Since $2t_{u} = u^{2} - 19u + 90$, the difference is $19(b-c)$, which accounts precisely for the factor $y^{b-c}$. Hence $\operatorname{ev}(\Theta)$ vanishes identically. Lastly, using $2t_{u} \equiv u^{2} + 14 \pmod{19}$, we obtain the stated weight congruence.
\end{proof}

\subsection{The doubling automorphism}

Let $\mathcal{T} \colon \mathcal{R} \to \mathcal{R}$ be defined on the variables by
\[
\mathcal{T}(Z_{x}) = (-1)^{x+1}Z_{2x} \qquad (x \in \mathbb{Z}).
\]
This agrees with \eqref{eq:Z-rules}: replacing $x$ by $19-x$ changes the sign $(-1)^{x+1}$ and replaces $Z_{2x}$ by $Z_{38-2x} = -Z_{2x}$, so the image is unchanged. Replacing $x$ by $-x$ maps $Z_{2x}$ to $-Z_{2x}$. Doubling permutes $(\mathbb{Z}/19)^{\times}/\{\pm1\}$ in the single $9$-cycle $1 \to 2 \to 4 \to 8 \to 3 \to 6 \to 7 \to 5 \to 9 \to 1$, so $\mathcal{T}$ permutes the variables up to sign and is a ring automorphism. For a monomial of degree $k$, the shift formula $t_{2u} = 4t_{u} - 135 + 19u \equiv 4t_{u}-2$ maps monomials of weight $w$ to weight $4w-2k$, respectively.

\begin{lemma}\label{lem:T-ideal}
$\mathcal{T}(\Theta(a,b,c,d)) = \Theta(2a,2b,2c,2d)$. Hence $\mathcal{T}(\mathcal{I}_{\Theta}) = \mathcal{I}_{\Theta}$.
\end{lemma}

\begin{proof}
Each monomial $Z_{x_{1}}Z_{x_{2}}Z_{x_{3}}Z_{x_{4}}$ of $\Theta$ has index sum $2a+2b$ or $2a+2c$, which is even. Thus, the sign product $\prod_{i}(-1)^{x_{i}+1}$ resulting from $\mathcal{T}$ is $+1$ for each term, while the index doubling maps $\Theta(a,b,c,d)$ directly to $\Theta(2a,2b,2c,2d)$.
\end{proof}

\begin{lemma}\label{lem:T-eigen}
$\mathcal{T}(C) = -\tfrac{1}{2}C$ and $\mathcal{T}(B) = \tfrac{1}{16}B$. Consequently, $\mathcal{T}(F) = 4F$ and $\mathcal{T}(F_{11}) = 4F_{11}$.
\end{lemma}

\begin{proof}
If $u$ is an integer, define $c(u) = (-1)^{u}2u\,Z_{u}$. By \eqref{eq:Z-rules}, since $c$ satisfies $c(-u) = c(u)$ and $c(u+19) = c(u)$, it follows that $c$ is a function on $(\mathbb{Z}/19)^{\times}/\{\pm1\}$ and $C = \sum c(u)$. Now,
\[
\mathcal{T}(c(u)) = (-1)^{u}2u \cdot (-1)^{u+1}Z_{2u} = -\tfrac{1}{2} \, c(2u),
\]
and doubling permutes the set.

Denoting $M(u,v) = Z_{u}Z_{v}Z_{u+v}Z_{v-u}$ for integers $u,v$ and $m(u,v) = (-1)^{u+v}uv(v^{2}-u^{2})M(u,v)$, one checks using \eqref{eq:Z-rules} that $m(v,u) = m(u,v)$, $m(-u,v) = m(u,v)$, and $m(u+19,v) = m(u,v)$. Thus, $m$ is a function on unordered pairs of distinct members of $(\mathbb{Z}/19)^{\times}/\{\pm1\}$, and $B$ is a sum over these pairs. Since the four coordinates of $M(u,v)$ sum to $u+3v$, we have $\mathcal{T}(M(u,v)) = (-1)^{u+v}M(2u,2v)$ and
\[
\mathcal{T}(m(u,v)) = uv(v^{2}-u^{2})M(2u,2v) = \tfrac{1}{16} \, m(2u,2v).
\]
Doubling permutes the pairs. Finally, combining the scaling factors, $\mathcal{T}(F) = \mathcal{T}(B \, C^{3}) = \frac{1}{16} \cdot \left(-\frac{1}{2}\right)^{3} B \, C^{3} = 4F$, and $\mathcal{T}(F_{11}) = 4F_{11}$ since $\mathcal{T}$ maps degree $7$ and weight $11$ to weight $44-14 \equiv 11$.
\end{proof}

\begin{remark}[The components in two orbits]\label{rmk:two-orbits}
Define $B^{(r)}$ for $0 \le r \le 18$ to be the sum of those terms of $B$ that have weight congruent to $r+30 \pmod{19}$, so that $\operatorname{ev}(B^{(r)}) = q^{r+30}B_{r}(y)$, as noted in the computation before Proposition~\ref{thm:bridge}. Since $\mathcal{T}$ maps elements of degree $4$ and weight $w$ to weight $4w - 8$, and $4(r+11) - 8 \equiv (4r+6)+11 \pmod{19}$, Lemma~\ref{lem:T-eigen} yields
\[
\mathcal{T}\bigl(B^{(r)}\bigr) = \tfrac{1}{16} \, B^{(4r+6)} \qquad (0 \le r \le 18).
\]
Thus, in the variables $Z_{a}$, the components $B_{r}$ in each orbit of Remark~\ref{rmk:doubling-rows} are obtained from any one of them through repeated index doubling, and Tables~\ref{tab:B} and~\ref{tab:Psi} are determined by their rows for $r=0$ and $r=2$. This reflects the action of the automorphism $\mathcal{T}$ on the components themselves.
\end{remark}

\begin{proposition}[Invariance of $S$]\label{prop:S}
Let $S_{Z} = Z_{1}Z_{7}Z_{8} - Z_{2}Z_{3}Z_{5} + Z_{4}Z_{6}Z_{9}$. Then $\mathcal{T}(S_{Z}) = S_{Z}$ in $\mathcal{R}$, and $\operatorname{ev}(S_{Z}) = q^{21}S$.
\end{proposition}

\begin{proof}
From the definition of $\mathcal{T}$ and the composition rules, we have
\begin{align*}
\mathcal{T}(Z_{1}Z_{7}Z_{8}) &= (+1)(+1)(-1) \, Z_{2}Z_{14}Z_{16} = -Z_{2}Z_{5}Z_{3}, \\
\mathcal{T}(Z_{2}Z_{3}Z_{5}) &= (-1)(+1)(+1) \, Z_{4}Z_{6}Z_{10} = -Z_{4}Z_{6}Z_{9}, \\
\mathcal{T}(Z_{4}Z_{6}Z_{9}) &= (-1)(-1)(+1) \, Z_{8}Z_{12}Z_{18} = Z_{8}Z_{7}Z_{1}.
\end{align*}
The monomials corresponding to these blocks have weights congruent to $21 \pmod{19}$. Specifically, evaluating under $\operatorname{ev}$ maps them directly to $q^{21}$ times $y P_{1}P_{7}P_{8}$, $y^{2}P_{2}P_{3}P_{5}$, and $P_{4}P_{6}P_{9}$, respectively, yielding the desired result.
\end{proof}

There are three blocks that additionally index three of the relations of Lemma~\ref{lem:theta}: for any block $\{b,c,d\}$, one has $b^{2}+c^{2}+d^{2} \equiv 0 \pmod{19}$ (since the sums are $114, 38,$ and $133$, which are multiples of $19$), and thus $\Theta(0,b,c,d)$ has weight $9$.

\subsection{Summary}

Collectively, the findings in this section can be summarized as follows. The value of the polynomial $F_{11}$ in Proposition~\ref{thm:scalar} evaluates to $q^{30}\mathcal{D}/P_{0}^{4}$ and satisfies $\mathcal{T}(F_{11}) = 4F_{11}$. Furthermore, $S_{Z}$ satisfies the condition
\[
\mathcal{T}(S_{Z}) = S_{Z},
\]
which evaluates to $q^{21}S$. Via O'Brien's identity \eqref{eq:obrien-factor}, we have
\[
\operatorname{ev}(F_{11}) \equiv \operatorname{ev}(S_{Z}) \cdot \frac{q^{9}(S^{2}-8y\mathcal{D}^{2})}{f(y)^{4}} \pmod{19}.
\]
The three products are permuted cyclically by the automorphism $\mathcal{T}$ acting on the cosets $\{1,7,8\}$, $\{2,3,5\}$, and $\{4,6,9\}$, as established in the proof of Proposition~\ref{prop:S}. 

In the span of these three products, the vectors fixed by $\mathcal{T}$ are exactly the multiples of $S_{Z}$. This does not show that $F_{11}$ is divisible by $S_{Z}$ in $\mathcal{R}$ modulo the ideal generated by the three-term relations; whether it is, is the first open question addressed in Section~\ref{sec:questions}.

\section{Questions}\label{sec:questions}

\begin{enumerate}
\item Is $F_{11}$ congruent to a quartic polynomial times $S_{Z}$ modulo the ideal $\mathcal{I}_{\Theta}$ generated by the three-term relations, and can this be shown without exhaustive computation? O'Brien proves the cubic relation \eqref{eq:obrien-cubic-Z} by restricting to a finite ansatz and matching coefficients. Can one give a proof directly in $\mathcal{R}$, and would the $\mathcal{T}$-covariance of Lemma~\ref{lem:T-eigen} account for this factorization algebraically? 

\item None of the components of $f(q)^{10}$ modulo $19$ are given by a single quartic product (Remark~\ref{rmk:contrast}). Can the two-term sums of Table~\ref{tab:B} be reduced to single products of the $P_{a}$ via relations among the $P_{a}$?

\item Let \eqref{eq:obrien-factor} be the relation given by Remark~\ref{rmk:weierstrass}, which can be viewed in terms of a Weierstrass model of $X_{0}(19)$ with nodal reduction. Does the cofactor $(S^{2}-8y\mathcal{D}^{2})/f(y)^{4}$ admit a natural geometric interpretation on that elliptic curve model?

\item By Lemma~\ref{lem:reduction} with $k=12$ and Proposition~\ref{prop:jacobi}, the superexceptional congruence $p_{7}(19n+9) \equiv 0$ \cite{KimingOlsson1992, Boylan2004} says exactly that the weight-$9$ part of $C^{4}$ vanishes on the theta functions; by Lemma~\ref{lem:T-eigen} it is an eigenvector of $\mathcal{T}$ with eigenvalue $1/16$. Is there an alternative conceptual proof of $p_{7}(19n+9) \equiv 0$, and can classical theta identities (such as specializations of Jacobi's five-term identity \cite[(1.2)]{Koornwinder2014}) be systematically applied here?
\end{enumerate}

\section*{Declarations}

\paragraph*{Funding.} This research received no external funding.

\paragraph*{Competing interests.} The author declares no competing interests.

\paragraph*{Ethics approval and consent to participate.} Not applicable.

\paragraph*{Data availability.} There are no data associated with this paper;
all results are proved in the text.

\paragraph*{Code availability.} Not applicable. No computer program is needed
for any proof in the paper.

\end{document}